\documentclass[12pt]{article}

\usepackage{amsmath,amssymb,amsthm,amscd,a4wide,upref}

\usepackage[enableskew]{youngtab}

\usepackage{hyperref}
\hypersetup{colorlinks,citecolor=blue,filecolor=black,linkcolor=blue,urlcolor=blue}
\usepackage{enumerate}
\usepackage{mathtools}
\usepackage{enumitem}

\usepackage{lmodern}     
\usepackage[T1]{fontenc}

\begin{document}

\newcommand{\End}{{\rm{End}\ts}}
\newcommand{\Hom}{{\rm{Hom}}}
\newcommand{\Mat}{{\rm{Mat}}}
\newcommand{\ch}{{\rm{ch}\ts}}
\newcommand{\sh}{{\rm{sh}}}
\newcommand{\chara}{{\rm{char}\ts}}
\newcommand{\diag}{ {\rm diag}}
\newcommand{\non}{\nonumber}
\newcommand{\wt}{\widetilde}
\newcommand{\wh}{\widehat}
\newcommand{\ot}{\otimes}
\newcommand{\la}{\lambda}
\newcommand{\La}{\Lambda}
\newcommand{\De}{\Delta}
\newcommand{\al}{\alpha}
\newcommand{\be}{\beta}
\newcommand{\ga}{\gamma}
\newcommand{\Ga}{\Gamma}
\newcommand{\ep}{\epsilon}
\newcommand{\ka}{\kappa}
\newcommand{\vk}{\varkappa}
\newcommand{\si}{\sigma}
\newcommand{\vp}{\varphi}
\newcommand{\de}{\delta}
\newcommand{\ze}{\zeta}
\newcommand{\om}{\omega}
\newcommand{\ee}{\epsilon^{}}
\newcommand{\su}{s^{}}
\newcommand{\hra}{\hookrightarrow}
\newcommand{\ve}{\varepsilon}
\newcommand{\ts}{\,}
\newcommand{\vac}{\mathbf{1}}
\newcommand{\vacr}{|\tss 0\rangle}
\newcommand{\vacl}{\langle 0\tss |}
\newcommand{\di}{\partial}
\newcommand{\qin}{q^{-1}}
\newcommand{\tss}{\hspace{1pt}}
\newcommand{\Sr}{ {\rm S}}
\newcommand{\U}{ {\rm U}}
\newcommand{\BL}{ {\overline L}}
\newcommand{\BE}{ {\overline E}}
\newcommand{\BP}{ {\overline P}}
\newcommand{\AAb}{\mathbb{A}\tss}
\newcommand{\CC}{\mathbb{C}\tss}
\newcommand{\KK}{\mathbb{K}\tss}
\newcommand{\QQ}{\mathbb{Q}\tss}
\newcommand{\SSb}{\mathbb{S}\tss}
\newcommand{\ZZ}{\mathbb{Z}\tss}
\newcommand{\X}{ {\rm X}}
\newcommand{\Y}{ {\rm Y}}
\newcommand{\Z}{{\rm Z}}
\newcommand{\Ac}{\mathcal{A}}
\newcommand{\achi}{\Ac_{\chi}}
\newcommand{\bachi}{\overline\Ac_{\chi}}
\newcommand{\Lc}{\mathcal{L}}
\newcommand{\Mc}{\mathcal{M}}
\newcommand{\Pc}{\mathcal{P}}
\newcommand{\Qc}{\mathcal{Q}}
\newcommand{\Tc}{\mathcal{T}}
\newcommand{\Sc}{\mathcal{S}}
\newcommand{\Bc}{\mathcal{B}}
\newcommand{\Dc}{\mathcal{D}}
\newcommand{\Ec}{\mathcal{E}}
\newcommand{\Fc}{\mathcal{F}}
\newcommand{\Hc}{\mathcal{H}}
\newcommand{\Uc}{\mathcal{U}}
\newcommand{\Vc}{\mathcal{V}}
\newcommand{\Wc}{\mathcal{W}}
\newcommand{\Yc}{\mathcal{Y}}
\newcommand{\Ar}{{\rm A}}
\newcommand{\Br}{{\rm B}}
\newcommand{\Ir}{{\rm I}}
\newcommand{\Fr}{{\rm F}}
\newcommand{\Jr}{{\rm J}}
\newcommand{\Or}{{\rm O}}
\newcommand{\GL}{{\rm GL}}
\newcommand{\Spr}{{\rm Sp}}
\newcommand{\Rr}{{\rm R}}
\newcommand{\Zr}{{\rm Z}}
\newcommand{\gl}{\mathfrak{gl}}
\newcommand{\middd}{{\rm mid}}
\newcommand{\ev}{{\rm ev}}
\newcommand{\Pf}{{\rm Pf}}
\newcommand{\Norm}{{\rm Norm\tss}}
\newcommand{\oa}{\mathfrak{o}}
\newcommand{\spa}{\mathfrak{sp}}
\newcommand{\osp}{\mathfrak{osp}}
\newcommand{\g}{\mathfrak{g}}
\newcommand{\h}{\mathfrak h}
\newcommand{\n}{\mathfrak n}
\newcommand{\z}{\mathfrak{z}}
\newcommand{\Zgot}{\mathfrak{Z}}
\newcommand{\p}{\mathfrak{p}}
\newcommand{\sll}{\mathfrak{sl}}
\newcommand{\agot}{\mathfrak{a}}
\newcommand{\qdet}{ {\rm qdet}\ts}
\newcommand{\Ber}{ {\rm Ber}\ts}
\newcommand{\HC}{ {\mathcal HC}}
\newcommand{\cdet}{ {\rm cdet}}
\newcommand{\tr}{ {\rm tr}}
\newcommand{\gr}{ {\rm gr}}
\newcommand{\str}{ {\rm str}}
\newcommand{\loc}{{\rm loc}}
\newcommand{\Gr}{{\rm G}}
\newcommand{\sgn}{ {\rm sgn}\ts}
\newcommand{\ba}{\bar{a}}
\newcommand{\bb}{\bar{b}}
\newcommand{\bi}{\bar{\imath}}
\newcommand{\bj}{\bar{\jmath}}
\newcommand{\bk}{\bar{k}}
\newcommand{\bl}{\bar{l}}
\newcommand{\hb}{\mathbf{h}}
\newcommand{\Sym}{\mathfrak S}
\newcommand{\fand}{\quad\text{and}\quad}
\newcommand{\Fand}{\qquad\text{and}\qquad}
\newcommand{\For}{\qquad\text{or}\qquad}
\newcommand{\OR}{\qquad\text{or}\qquad}

\renewcommand{\theequation}{\arabic{section}.\arabic{equation}}

\newtheorem{thm}{Theorem}[section]
\newtheorem{lem}[thm]{Lemma}
\newtheorem{prop}[thm]{Proposition}
\newtheorem{cor}[thm]{Corollary}
\newtheorem{conj}[thm]{Conjecture}
\newtheorem*{mthm}{Main Theorem}
\newtheorem*{mthma}{Theorem A}
\newtheorem*{mthmb}{Theorem B}

\theoremstyle{definition}
\newtheorem{defin}[thm]{Definition}

\theoremstyle{remark}
\newtheorem{remark}[thm]{Remark}
\newtheorem{example}[thm]{Example}

\newcommand{\bth}{\begin{thm}}
\renewcommand{\eth}{\end{thm}}
\newcommand{\bpr}{\begin{prop}}
\newcommand{\epr}{\end{prop}}
\newcommand{\ble}{\begin{lem}}
\newcommand{\ele}{\end{lem}}
\newcommand{\bco}{\begin{cor}}
\newcommand{\eco}{\end{cor}}
\newcommand{\bde}{\begin{defin}}
\newcommand{\ede}{\end{defin}}
\newcommand{\bex}{\begin{example}}
\newcommand{\eex}{\end{example}}
\newcommand{\bre}{\begin{remark}}
\newcommand{\ere}{\end{remark}}
\newcommand{\bcj}{\begin{conj}}
\newcommand{\ecj}{\end{conj}}

\newcommand{\bal}{\begin{aligned}}
\newcommand{\eal}{\end{aligned}}
\newcommand{\beq}{\begin{equation}}
\newcommand{\eeq}{\end{equation}}
\newcommand{\ben}{\begin{equation*}}
\newcommand{\een}{\end{equation*}}

\newcommand{\bpf}{\begin{proof}}
\newcommand{\epf}{\end{proof}}

\def\beql#1{\begin{equation}\label{#1}}

\title{\Large\bf Quantum Sugawara operators in type $A$}

\author{{Naihuan Jing,\ \   Ming Liu\ \   and\ \   Alexander Molev}}

\date{} 
\maketitle

\vspace{4 mm}

\begin{abstract}
We construct Sugawara operators for the quantum affine algebra of type $A$ in an explicit form.
The operators are associated with primitive idempotents of the Hecke algebra and
parameterized by Young diagrams. This generalizes a previous construction (2016)
where one-column diagrams were considered. We calculate the Harish-Chandra
images of the Sugawara operators and identify them
with the eigenvalues of the operators acting in the $q$-deformed Wakimoto modules.

%

\end{abstract}


\vspace{5 mm}

%


\section{Introduction}
\label{sec:int}
\setcounter{equation}{0}

The {\em Sugawara operators} associated with a simple Lie algebra $\g$
are elements of the center of a completed universal enveloping algebra
$\wt\U(\wh\g)$ of the affine Lie algebra $\wh\g$ at the critical level.
A detailed account of Sugawara operators based on the vertex algebra theory
along with their connections with the geometry of opers is given in the book \cite{f:lc}; see also
\cite{m:so} for explicit constructions in the classical types.

The $q$-versions of the operators originate in \cite{rs:ce}; these are
elements of the center $\Zr_q(\wh\g)$ of the completion $\wt\U_q(\wh\g)$
of the quantum affine algebra at the critical level and we call them
the {\em quantum Sugawara operators} in this paper. According to \cite{rs:ce},
a formal Laurent series $\ell_V(z)$ with coefficients in $\Zr_q(\wh\g)$
can be associated with every finite-dimensional representation $V$
of the quantum affine algebra $\U_q(\wh\g)$. By the results of \cite{de:cq},
the map $V\mapsto \ell_V(z)$ possesses properties of a homomorphism from
the Grothendieck ring $\text{\rm Rep}\ts \U_q(\wh\g)$ to formal series in $z$ and
the coefficients of $\ell_V(z)$
generate all singular vectors in Verma modules. Furthermore, the Harish-Chandra images
of the quantum Sugawara operators essentially coincide with the $q$-{\em characters}
of finite-dimensional representations of $\U_q(\wh\g)$ and
can also be interpreted as elements of the $q$-{\em deformed
classical $\Wc$-algebras}; see \cite{fr:qc}.

A key starting point in the above mentioned constructions of the quantum Sugawara
operators is the existence of the {\em universal
$R$-matrix} in the Drinfeld presentation of $\U_q(\wh\g)$ \cite{d:qg}. The series $\ell_V(z)$
is defined via the $L$-{\em operators} constructed from the universal
$R$-matrix and its properties imply the centrality property of $\ell_V(z)$.

In this paper we take a different viewpoint and {\em define}
the quantum affine algebra $\U_q(\wh\g)$ in type $A$ via its $RLL$ presentation.
Given an arbitrary Young
diagram $\la$ containing at most $n$ rows, we
produce explicit formulas for
Laurent series $S_{\La}(z)$ in $z$, whose coefficients are
quantum Sugawara operators (Theorem~\ref{thm:centr}). They are elements of
the center of the completed universal enveloping algebra $\wt\U_q(\wh\gl_n)$
at the critical level. In the particular case
of one-column diagrams $\la$ these operators were already found in \cite{fjmr:hs}.

Our arguments are based on the quantum Schur--Weyl duality \cite{j:qu} between the
Hecke algebra $\Hc_m$ and the quantized enveloping algebra $\U_q(\gl_n)$.
Given a Young diagram $\la$ with $m$ boxes, we choose a standard $\la$-tableau $\La$ and
consider the corresponding primitive idempotent $e^{\la}_{\La}\in \Hc_m$. By employing the commuting
actions of $\Hc_m$ and $\U_q(\gl_n)$ on the vector space $(\CC^n)^{\ot m}$, we define
the series $S_{\La}(z)$ as the weighted trace of the image of
$e^{\la}_{\La}$ in the representation,
multiplied by a product of $L$-operators; see formula \eqref{lkz} below. The series $S_{\La}(z)$
turns out not to depend on the choice of $\La$.
The images of the coefficients of $S_{\La}(z)$ in the vacuum module at the critical level are
invariants of the module (Corollary~\ref{cor:inv}).

In the proof that the coefficients of $S_{\La}(z)$ belong to the center of $\wt\U_q(\wh\gl_n)$
at the critical level, we generalize the $R$-matrix calculations of \cite{fjmr:hs} and
use the {\em fusion procedure} for the Hecke algebra \cite{c:ni}, \cite{n:mh}; see also \cite{imo:ih}.
Although the part of the proof involving the Hecke algebra relies on suitable $q$-versions of
the symmetric group arguments originated in \cite{o:qi}, some significant new ingredients turned out
to be necessary for the evaluation of the Harish-Chandra images; cf.
the Yangian versions considered in \cite{jkmy:cq} and \cite[Ch.~10]{m:so}.

We will show that, upon a suitable identification
of the parameters, the Harish-Chandra images of the series $S_{\La}(z)$ coincide with the
$q$-characters of the evaluation modules $L(\la^{\circ})$ over $\U_q(\wh\gl_n)$ associated with
$\la$ (see Remark~\ref{rem:qchar}).
Moreover, the coefficients of the series $S_{\La}(z)$ are known to act as multiplication
by scalars in the $q$-deformed Wakimoto modules
over $\U_q(\wh\gl_n)$ at the critical level \cite{aos:fb}. We show that
the eigenvalues essentially coincide with the Harish-Chandra images
of the coefficients. This extends the
corresponding results of \cite{fjmr:hs} to arbitrary Young diagrams.

\section{Schur elements and fusion procedure for Hecke algebras}
\label{sec:se}
\setcounter{equation}{0}

Throughout the paper we assume that $q$ is a fixed nonzero complex number which is not a root
of unity. We let $\Hc_m$ denote
the {\em Hecke algebra} over $\CC$ which is generated by elements
$T_1,\dots,T_{m-1}$ subject to the relations
\ben
\bal
(T_i-q)(T_i+q^{-1})=0,\\
T_{i}T_{i+1}T_{i}=T_{i+1}T_{i}T_{i+1},\\
T_iT_j=T_jT_i \quad\text{for\ \  $|i-j|>1$}.
\eal
\een
The condition on $q$ ensures that the algebra $\Hc_m$ is semisimple and we will recall some
basic facts about its representations; see e.g. \cite{gp:cf}, \cite{m:rt}.

For $i=1, \dots, m-1$
we let $\si_{i}=(i, i+1)$ be the adjacent transpositions in the symmetric group $\Sym_{m}$.
Choose a reduced decomposition
$\si=\si_{i_{1}} \dots \si_{i_{l}}$ of any element
$\si \in \Sym_{m}$ and set $T_{\si}=T_{i_{1}} \dots T_{i_{l}}$. This element of
the algebra $\Hc_{m}$ does not depend on
the choice of reduced decomposition of $\si$. The element of maximal length
in $\Sym_{m}$ will be denoted by $\si_{0}$ and we set $T_{0}=T_{\si_{0}}$.
The elements $T_{\si}$ parameterized by $\si\in\Sym_m$ form a basis of $\Hc_{m}$.
Set $T_{\si}^{*}=T_{\si^{-1}}$ and extend this star operation
to $\Hc_m$ by linearity. It defines an involutive
anti-automorphism of $\Hc_m$.

\paragraph{Young basis and matrix units.}
We will use some standard notation on partitions. Any partition
$\la=(\la_{1}, \dots, \la_{l})$ will be identified with the corresponding Young diagram
which is a left-justified array of rows of unit boxes such that the top row
contains $\la_1$ boxes, the next row contains $\la_2$ boxes, etc.
We will need the conjugate partition $\la^{\prime}=(\la_{1}^{\prime}, \dots, \la_{r}^{\prime})$
so that $\la_{j}^{\prime}$ is the number of boxes in the $j$-th column of
the Young diagram $\la$. If $\alpha=(i, j)$
is a box of $\la,$ then the corresponding hook length
is defined as $h(\al)=\la_{i}+\la_{j}^{\prime}-i-j+1$ and
the content is $c(\al)=j-i$.

The irreducible representations of $\Hc_{m}$ over $\CC$ are parameterized by partitions of $m$.
Given such a partition $\la \vdash m$, we let $V_{\la}$ denote the
corresponding irreducible representation and let
\beql{phila}
\vp_{\la}:\Hc_{m}\to \End V_{\la}
\eeq
be the associated homomorphism.
The vector space $V_{\la}$ is equipped with an $\Hc_{m}$-invariant
inner product $\langle\ts\tss,\ts\rangle=\langle\ts\tss,\ts\rangle^{}_{V_{\la}}$ so that
$\langle hv_1,v_2\rangle=\langle v_1,h^*v_2\rangle$ for $h\in \Hc_m$ and $v_1,v_2\in  V_{\la}$.

A tableau $\Uc$
of shape $\la\vdash m$
(or a $\la$-tableau)
is obtained by
filling in the boxes of the diagram with the numbers
in a given set $\{1,\dots,n\}$. The tableau is called {\em semistandard}
if the entries weakly increase along each row from left to right
and strictly increase in each column from top to bottom.
We write $\sh(\Uc)=\la$ if the shape of $\Uc$ is $\la$.
A $\la$-tableau with entries in
$\{1,\dots,m\}$ which are filled in the boxes bijectively
is called {\em standard}
if its entries strictly increase along the rows and down the columns.

The {\em orthonormal Young basis} $\{v_{\La}\}$
of $V_{\la}$ is parameterized by the set of {\em standard $\la$-tableaux} $\La$.
We will write $f_{\la}=\dim V_{\la}$ for the number of such tableaux.
Given a standard $\la$-tableau $\La$, we let $c_{k}(\La)$ denote the content $j-i$
of the box $(i,j)$ of $\la$ occupied by $k$ in $\La$. We also set $d_k(\La)=c_{k+1}(\La)-c_{k}(\La)$
and use the standard notation
\ben
[n]_{q}=\frac{q^{n}-q^{-n}}{q-q^{-1}}\Fand [n]_{q}!=[1]_{q}\cdots [n]_{q}.
\een
The action
of the generators of $\Hc_{m}$ in the orthonormal Young basis of $V_{\la}$
can be described explicitly by using the {\em seminormal form} of \cite{dj:bi} and \cite{h:rh}.
Namely, by normalizing the basis vectors,
for any
$k \in\{1, \dots, m-1\}$ we get
\beql{eq:Young basis}
T_{k} \ts v^{}_{\La}=\frac{q^{d_k(\La)}}{[d_k(\La)]_q} v^{}_{\La}+\sqrt{1-\frac{1}{[d_k(\La)]_q^2}}
\ts v^{}_{\si_{k} \La},
\eeq
where the tableau $\si_{k} \La$ is obtained from $\La$ by swapping the entries $k$ and $k+1$, and
we assume that $v_{\si_{k} \La}=0$ if the tableau $\si_{k} \La$ is not standard.
We will denote the character of the $\Hc_m$-module
$V_{\la}$ by $\chi_{\la}$.

For any skew diagram $\theta$ with $m$ boxes,
the same formulas \eqref{eq:Young basis} define a representation of the Hecke algebra $\Hc_{m}$
on the linear span $V_{\theta}$ of vectors $v_{\La}$ parameterized by standard tableaux $\La$
of shape $\theta$. We will write $\chi_{\theta}$ for the character of $V_{\theta}$.

The {\em symmetrizing trace} on the Hecke algebra is the linear map
\ben
\tau:\Hc_{m}\to\CC,\qquad \tau(T_{\si})=\de_{\si,e},
\een
where $e$ denotes the identity element of $\Sym_m$. This induces a {\em symmetric algebra}
structure on $\Hc_{m}$ associated with the bilinear form
\ben
\Hc_{m}\ot\Hc_m\to\CC,\qquad h_1\ot h_2\mapsto\tau(h_1h_2),
\een
which is symmetric and non-degenerate.  The dual basis of
$\{T_{\si}\}_{\si \in \Sym_m}$ with respect to the form is
$\{ T_{\si^{-1}}\}_{\si \in \Sym_m}$ \cite[Prop.~4.45]{m:rt}.

For any $u\in\End V_{\la}$, define $I(u)\in \End V_{\la}$ by
\ben
I(u)=\sum_{\si\in \Sym_m}\vp_\la(T_{\si})\ts u\ts\vp_\la(T_{\si^{-1}}).
\een
It is straightforward to verify that the operator
$I(u)$ belongs to ${\rm{End}}_{\Hc_m}V_{\la}$ and therefore must be proportional to
the identity operator by Schur's Lemma. In fact,
\beql{iuschur}
I(u)=c_{\la}\tss\tr^{}_{V_{\la}}(u)\tss{\rm id}^{}_{V_{\la}},
\eeq
where the constant factor $c_{\la}$ is the {\em Schur element} given by
the Steinberg formula
\ben
c_{\la}=\prod_{\al\in\la} q^{c(\al)}\tss [h(\al)]_q;
\een
see e.g. \cite[Thm~4.64]{m:rt} for the calculation. We get the
following formula for the symmetrizing trace
\beql{tauchi}
\tau=\sum_{\la\vdash m}\frac{\chi_{\la}}{c_{\la}}.
\eeq

The Hecke algebra $\Hc_m$ is semisimple and isomorphic to the direct sum of matrix algebras
\beql{isomhecke}
\Hc_m \cong \bigoplus_{\la \vdash m} \Mat_{f_{\la}}(\CC).
\eeq
Via this isomorphism, the matrix units $e^{\la}_{\La,\Ga} \in \Mat_{f_{\la}}(\CC)$
can be identified with elements of $\Hc_m$ by
\beql{mat-units}
e^{\la}_{\La,\Ga}=\frac{1}{c_{\la}}\sum_{\si\in \Sym_m}
\langle T_{\si^{-1}}v^{}_{\La},v^{}_{\Ga}\rangle^{}_{V_{\la}}T_{\si}.
\eeq
We also point out the character orthogonality relations
\beql{orthcha}
\sum_{\si\in \Sym_m}\chi_{\la}(T_{\si})\chi_{\mu}(T_{\si^{-1}})=\de_{\la,\mu}c_{\la}f_{\la}
\eeq
which hold for any partitions $\la,\mu\vdash m$.

\paragraph{Fusion procedure.}
The diagonal matrix units $e^{\la}_{\La}=e^{\la}_{\La \La}$ are primitive idempotents of $\Hc_m$.
They will play a special role in what follows so we will produce a few formulas
for their calculation. Note their immediate properties:
\ben
e^{\la}_{\La} e^{\la}_{\Ga}=0 \quad \text { if } \quad \La \neq \Ga\Fand
(e^{\la}_{\La})^{2}=e^{\la}_{\La}
\een
together with the decomposition of the identity in $\Hc_m$
\ben
1=\sum_{\la \vdash m} \sum_{\sh(\La)=\la} e^{\la}_{\La}.
\een

\ble\label{lem:et}
Let $\La$ be a standard $\la$-tableau and let  $k \in\{1, \dots, m-1\}$.
Then
\ben
e^{\la}_{\La}\left(T_{k}-\frac{q^{d_k(\La)}}{[d_k(\La)]_q}\right)=\sqrt{1-\frac{1}{[d_k(\La)]_q^2}}
\ts e^{\la}_{\La,\si_{k} \La}
\een
and
\ben
\left(T_{k}+\frac{q^{-d_k(\La)}}{[d_k(\La)]_q}\right)e^{\la}_{\si_{k} \La}
=\sqrt{1-\frac{1}{[d_k(\La)]_q^2}} \ts e^{\la}_{\La,\si_{k} \La},
\een
assuming that $e^{\la}_{\La,\si_{k} \La}=0$ if the tableau $\si_{k} \La$ is not standard.
\ele

\bpf
Both relations are verified in the same way, so we will only consider the first one. Write
\ben
e^{\la}_{\La}T_{k}=\sum_{\si\in \Sym_m}a(\si)\ts T_{\si},\qquad a(\si)\in\CC.
\een
By applying \eqref{tauchi} we get
\ben
a(\si)=\tau(e^{\la}_{\La}T_kT_{\si^{-1}})
=\sum_{\mu\vdash m}\frac{\chi_{\mu}(e^{\la}_{\La}T_kT_{\si^{-1}})}{c_{\mu}}
=\frac{1}{c_{\la}}\chi_{\la}(e^{\la}_{\La}T_kT_{\si^{-1}})
=\frac{1}{c_{\la}}\chi_{\la}(T_{\si^{-1}}e^{\la}_{\La}T_k).
\een
Since
\ben
\chi_{\la}(u)=\sum_{\Ga,\ts\sh(\Ga)=\la}\langle u\ts
v^{}_{\Ga},v^{}_{\Ga}\rangle,
\een
we can use \eqref{eq:Young basis} and the matrix unit formula \eqref{mat-units}
to bring this to the form
\ben
\bal
a(\si)&=\frac{1}{c_{\la}}\sum_{\sh(\Ga)=\la}\frac{q^{d_k(\Ga)}}{[d_k(\Ga)]_q}
\big\langle T_{\si^{-1}}e^{\la}_{\La}v^{}_{\Ga},v^{}_{\Ga}\big\rangle
+\frac{1}{c_{\la}}\sum_{\sh(\Ga)=\la}\sqrt{1-\frac{1}{[d_k(\Ga)]_q^2}}\ts
\big\langle T_{\si^{-1}}e^{\la}_{\La}v^{}_{\si_k\Ga},v^{}_{\Ga}\big\rangle\\[0.4em]
&=\frac{1}{c_{\la}}\frac{q^{d_k(\La)}}{[d_k(\La)]_q}\big
\langle T_{\si^{-1}}v^{}_{\La},v^{}_{\La}\big\rangle
+\frac{1}{c_{\la}}\sqrt{1-\frac{1}{[d_k(\La)]_q^2}}\ts
\big\langle T_{\si^{-1}}v^{}_{\La},v^{}_{\si_k\La}\big\rangle,
\eal
\een
thus yielding the required formula.
\epf

The primitive idempotents  $e^{\la}_{\La}$ can be expressed explicitly in terms
of the pairwise commuting {\em Jucys--Murphy elements} of $\Hc_{m}$ which are defined by $y_1=1$ and
\ben
y_k=1+(q-\qin)\ts
\big(T_{(1,k)}+T_{(2,k)}+\dots+T_{(k-1,k)}\big),\qquad k=2,\dots,m,
\een
where the elements $T_{(i,j)}\in\Hc_m$ are
associated with the transpositions $(i,j)\in\Sym_m$; see
\cite{c:ni}, \cite{dj:bi}.
Let $m\geqslant 2$ and let $\la$ be a partition of $m$.
Fix a standard $\la$-tableau $\La$ and denote by
$\Ga$ the standard tableau
obtained from $\La$ by removing the box $\al$ occupied by $m$.
Then the shape of $\Ga$ is a diagram which we denote by $\mu$.
We let $c$ denote the content of the box $\al$ and
let $u$ be a complex variable. We have the recurrence relations
\ben
e^{\la}_{\La}=e^{\mu}_{\Ga}\ts \frac{u-q^{2c}\ }{u-y_m}\ts\Big|^{}_{u=q^{2c}}
\Fand
e^{\la}_{\La}=e^{\mu}_{\Ga} \ts\frac{(y_{m}-q^{2a_{1}})
\dots(y_{m}-q^{2a_{l}})}{(q^{2c}-q^{2a_{1}}) \dots(q^{2c}-q^{2a_{l}})},
\een
where $a_{1}, \dots, a_{l}$ are the contents of all
addable boxes of $\mu$ except for $\al$, while $c$ is the content of the latter.

Another way to calculate the primitive idempotents is based on the {\em fusion procedure}
for the Hecke algebra $\Hc_m$, originated in \cite{c:ni} with detailed proofs given in \cite{n:mh}.
It will be sufficient for us to use its simple version obtained in \cite{imo:ih} which we recall
below. For each $k=1,\dots, m-1$, introduce the $\Hc_{m}$-valued rational functions in two variables
$x, y$ by
\ben
T_{k}(x,y)=T_k+\frac{q-q^{-1}}{x^{-1}y-1},
\een
which satisfy the relations
\ben
T_{k}(x,y)\ts T_{k+1}(x,z)\ts T_{k}(y,z)=T_{k+1}(y,z)\ts T_{k}(x,z)\ts T_{k+1}(x,y)
\een
and
\ben
T_{k}(x, y)\ts T_{k}(y, x)=1-\frac{\left(q-q^{-1}\right)^{2} x y}{(x-y)^{2}}.
\een

Let $\La$ be a standard tableau of shape $\la$. We will keep the notation $c_{k}(\La)$
for the content $j-i$ of the box $(i,j)$ of $\la$ occupied by $k$ in $\La$.
Equip the set of all pairs $(i, j)$
with $1 \leqslant i<j \leqslant m$ with the following ordering. The pair $(i, j)$ precedes
$(i^{\tss\prime}, j^{\tss\prime})$ if $j<j^{\tss\prime}$, or if $j=j^{\tss\prime}$
but $i<i^{\tss\prime}$. Set
\beql{F}
T_{\La}(z_{1}, \dots, z_{m})=
\overrightarrow{\prod\limits_{(i, j)}}\ts
T_{j-i}\big(q^{2 c_{i}(\La)} z_{i}, q^{2 c_{j}(\La)} z_{j}\big)
\eeq
with the ordered product taken over the set of pairs. This is a rational function
in variables $z_{1}, \dots, z_{m}$ taking values
in $\Hc_{m}$. According to \cite[Thm~3.3]{imo:ih},
the primitive idempotent $e^{\la}_{\La}$ can be obtained by
the consecutive evaluations
\beql{theor}
e^{\la}_{\La} = \frac{1}{c_{\la'}}\ts T_{\La}(z_{1}, \dots, z_{m})\ts T_0^{-1}\ts
\big|_{z_1=1}\big|_{z_2=1} \dots \big|_{z_m=1},
\eeq
which are all well-defined; the rational functions are regular at
the evaluation points at each step, and $c_{\la'}$ is the Schur element associated with $\la'$.

\section{Formulas for Sugawara operators}
\label{sec:fs}
\setcounter{equation}{0}

As before, we regard $q$ as a nonzero complex number which is not
a root of unity. We will take the
$RLL$ presentation of the quantum affine algebra $\U_q(\wh\gl_n)$ introduced in \cite{rs:ce}
as its {\em definition}. We let $e_{ij}\in \End\CC^n$ denote the standard matrix units and
consider the $R$-matrices
\begin{multline}
\overline R(x)=\sum_{i}e_{ii}\ot e_{ii}\\[0.4em]
{}+{}\frac{1-x}{q-\qin x}\ts\sum_{i\ne j}e_{ii}\ot e_{jj}
+\frac{(q-\qin)\ts x}{q-\qin x}\ts\sum_{i> j}e_{ij}\ot
e_{ji}+ \frac{q-\qin}{q-\qin x}\ts\sum_{i< j}e_{ij}\ot e_{ji}
\label{rbar}
\end{multline}
and
\beql{rf}
R(x)=f(x)\tss \overline R(x),
\eeq
where
\ben
f(x)=1+\sum_{k=1}^{\infty}f_kx^k,\qquad f_k=f_k(q),
\een
is a formal power series in $x$ whose coefficients $f_k$ are
uniquely determined
by the relation
\ben
f(xq^{2n})=f(x)\ts\frac{(1-xq^2)\tss(1-xq^{2n-2})}{(1-x)\tss(1-xq^{2n})}.
\een
The {\em quantum affine algebra $\U_q(\wh\gl_n)$}
is generated by elements
\ben
l^+_{ij}[-r],\qquad l^-_{ij}[r]\qquad\text{with}\quad 1\leqslant i,j\leqslant n,\qquad r=0,1,\dots,
\een
and the invertible central element $q^c$,
subject to the defining relations
\begin{align}
l^+_{ji}[0]&=l^-_{ij}[0]=0\qquad&&\text{for}\qquad 1\leqslant i<j\leqslant n,
\non\\
l^+_{ii}[0]\ts l^-_{ii}[0]&=l^-_{ii}[0]\ts l^+_{ii}[0]=1\qquad&&\text{for}\qquad i=1,\dots,n,
\non
\end{align}
and
\begin{align}
R(u/v)L_1^{\pm}(u)L_2^{\pm}(v)&=L_2^{\pm}(v)L_1^{\pm}(u)R(u/v),
\label{RLL}\\[0.2em]
R(uq^{-c}/v)L_1^{+}(u)L_2^{-}(v)&=L_2^{-}(v)L_1^{+}(u)R(uq^{c}/v).
\label{RLLpm}
\end{align}
In the last two relations we consider the matrices $L^{\pm}(u)=\big[\tss l^{\pm}_{ij}(u)\big]$,
whose entries are formal power series in $u$ and $u^{-1}$,
\beql{serlpm}
l^{+}_{ij}(u)=\sum_{r=0}^{\infty}l^{+}_{ij}[-r]\tss u^r,\qquad
l^{-}_{ij}(u)=\sum_{r=0}^{\infty}l^{-}_{ij}[r]\tss u^{-r}.
\eeq
Here and below we regard the matrices as elements
\ben
L^{\pm}(u)=\sum_{i,j=1}^n e_{ij}\ot l^{\pm}_{ij}(u)\in\End\CC^n\ot\U_q(\wh\gl_n)[[u^{\pm1}]]
\een
and use a subscript to indicate a copy of the matrix in the multiple
tensor product algebra
\beql{multtpr}
\underbrace{\End\CC^n\ot\dots\ot\End\CC^n}_m\ot\U_q(\wh\gl_n)[[u^{\pm1}]]
\eeq
so that
\ben
L^{\pm}_a(u)=\sum_{i,j=1}^n 1^{\ot (a-1)}\ot e_{ij}\ot 1^{\ot (m-a)}\ot l^{\pm}_{ij}(u).
\een
In particular, we take $m=2$ for the defining relations \eqref{RLL} and \eqref{RLLpm}.

This notation for elements of algebras of the form \eqref{multtpr} will be extended
as follows. For an element
\ben
C=\sum_{i,j,r,s=1}^n c^{}_{ijrs}\ts e_{ij}\ot e_{rs}\in
\End \CC^n\ot\End \CC^n
\een
and any two indices $a,b\in\{1,\dots,m\}$ such that $a\ne b$,
we denote by $C_{a\tss b}$ the element of the algebra $(\End\CC^n)^{\ot m}$ with $m\geqslant 2$
given by
\beql{cab}
C_{a\tss b}=\sum_{i,j,r,s=1}^n c^{}_{ijrs}\ts (e_{ij})_a\tss (e_{rs})_b,
\qquad
(e_{ij})_a=1^{\ot(a-1)}\ot e_{ij}\ot 1^{\ot(m-a)}.
\eeq
We regard the matrix transposition as the linear map
\ben
t:\End\CC^n\to\End\CC^n,\qquad e_{ij}\mapsto e_{ji}.
\een
For any $a\in\{1,\dots,m\}$ we will denote by $t_a$ the corresponding
partial transposition on the algebra \eqref{multtpr} which acts as $t$ on the
$a$-th copy of $\End \CC^n$ and as the identity map on all the other tensor factors.

The $R$-matrix \eqref{rf} satisfies the {\em crossing symmetry relations} \cite{fri:qa}:
\beql{cs}
\big(R_{12}(x)^{-1}\big)^{t_2} D_2 R_{12}(xq^{2n})^{t_2}=D_2
\Fand
R_{12}(xq^{2n})^{t_1}\tss D_1\big(R_{12}(x)^{-1}\big)^{t_1}=D_1,
\eeq
where $D$ denotes the diagonal $n\times n$ matrix
\beql{d}
D=\diag\big[q^{n-1}, q^{n-3},\dots, q^{-n+1}\big]
\eeq
with the meaning of subscripts as in \eqref{cab}.

Denote by $\U_q(\wh\gl_n)_{\text{\rm cri}}$ the quantum affine algebra
{\em at the critical level $c=-n$}, which is
the quotient of $\U_q(\wh\gl_n)$ by the relation $q^c=q^{-n}$.
Its completion $\wt\U_q(\wh\gl_n)_{\text{\rm cri}}$ is defined as
the inverse limit
\beql{compl}
\wt\U_q(\wh\gl_n)_{\text{\rm cri}}=\lim_{\longleftarrow}
\U_q(\wh\gl_n)_{\text{\rm cri}}/J_p, \qquad p>0,
\eeq
where $J_p$ denotes the left ideal of $\U_q(\wh\gl_n)_{\text{\rm cri}}$ generated by all elements
$l^{-}_{ij}[r]$ with $r\geqslant p$.
Elements of the center $\Zr_q(\wh\gl_n)$ of $\wt\U_q(\wh\gl_n)_{\text{\rm cri}}$
are known as ({\em quantum}) {\em Sugawara operators}.

The $R$-matrix $\overline R(x)$ can be written as
\ben
\overline R(x)=\frac{1-x}{q-q^{-1}x}\ts\Big(R+\frac{q-q^{-1}}{x^{-1}-1}P\Big),
\een
where
\beql{rfin}
R=q\sum_{i=1}^{n}e_{ii}\ot e_{ii}+\sum_{i\neq j}e_{ii}\ot e_{jj}+(q-q^{-1})\sum_{i<j}e_{ij}\ot e_{ji},
\eeq
and
\ben
P=\sum\limits_{i,j=1}^{n}e_{ij}\ot e_{ji}.
\een
Setting $\check{R}=PR$, we get
\ben
\check{R}_{k}\check{R}_{k+1}\check{R}_{k}=\check{R}_{k+1}\check{R}_{k}\check{R}_{k+1}\Fand
(\check{R}_{k}-q)(\check{R}_{k}+q^{-1})=0,
\een
where $\check{R}_{k}=P_{k,k+1}R_{k,k+1}$.
Hence we obtain a representation of the Hecke algebra $\Hc_{m}$ on the tensor product space
$(\CC^{n})^{\ot m}$
defined by
\beql{haact}
T_{k}\mapsto \check{R}_{k},\qquad k=1,\dots,m-1.
\eeq
We will write $\check{R}_{0}$ for the image of the element $T_0\in \Hc_{m}$.
Furthermore, we find that under this action of the Hecke algebra,
\ben
T_{k}(x,y)\mapsto \check{R}_{k}(x/y)\qquad\text{with}\quad
\check{R}_{k}(z)=\check{R}_k+\frac{q-q^{-1}}{z^{-1}-1}.
\een
Given a standard tableau $\La$ of shape $\la\vdash m$,
we will let $\check{R}_{\La}(z_{1}, \dots, z_{m})$ denote the image of
the product \eqref{F} and
denote by $\Ec^{\la}_{\La}$ (or just $\Ec_{\La}$)  the image of the primitive idempotent $e^{\la}_{\La}$
under the action \eqref{haact}. Clearly, $\Ec_{\La}^2=\Ec_{\La}$, and the fusion formula
\eqref{theor} implies the relation
\beql{fusr}
\Ec_{\La} = \frac{1}{c_{\la'}}\ts \check{R}_{\La}(z_{1}, \dots, z_{m})\ts \check{R}_{0}^{-1}\ts
\big|_{z_1=1}\big|_{z_2=1} \dots \big|_{z_m=1}.
\eeq
The operator $\Ec_{\La}$ is zero unless the Young diagram $\la$ contains at most $n$ rows.
In what follows, we will assume that $\la$ satisfies this condition.
Introduce the Laurent series $S_{\La}(z)$ in $z$ by
\beql{lkz}
\bal
S_{\La}(z)=\tr^{}_{1,\dots,m}\ts &L^+_1(zq^{-2c_1(\La)})\dots L^+_m(zq^{-2c_{m}(\La)})\\
&\times L^-_m(zq^{-n-2c_{m}(\La)})^{-1}\dots L^-_1(zq^{-n-2c_1(\La)})^{-1}D_1\dots D_m \Ec_{\La} ,
\eal
\eeq
where $D$ is the diagonal matrix \eqref{d} and
the trace is taken over all $m$ copies of $\End\CC^n$ in \eqref{multtpr}.
All coefficients of the series $S_{\La}(z)$
are elements of the algebra $\wt\U_q(\wh\gl_n)_{\text{\rm cri}}$.

The following is our main result which provides explicit formulas
for quantum Sugawara operators.

\bth\label{thm:centr}
The coefficients of the series $S_{\La}(z)$ belong to the center $\Zr_q(\wh\gl_n)$ of the completed quantum
affine algebra at the critical level $\wt\U_q(\wh\gl_n)_{\text{\rm cri}}$. Moreover,
$S_{\La}(z)$ does not depend on the standard $\la$-tableau $\La$
and only depends on the Young diagram $\la$.
\eth

The proof of Theorem~\ref{thm:centr} will be given in the rest of this section.
Set
\ben
L^{\pm}_{\La}(z)=L^\pm_1(zq^{-2c_1(\La)})\dots L^\pm_m(zq^{-2c_{m}(\La)})
\een
and $L_{\La}(z)=L^{+}_{\La}(z)L^{-}_{\La}(zq^{-n})^{-1}$. Then
$S_{\La}(z)$ can be written as
\ben
S_{\La}(z)=\tr^{}_{1,\dots,m}\ts L_{\La}(z)D_1\dots D_m \Ec_{\La}.
\een

\ble\label{lem:LE}
We have the relations
\ben
L^{+}_{\La}(z)\tss\Ec_{\La}=\Ec_{\La}\tss L^{+}_{\La}(z)\tss\Ec_{\La},
\een
\ben
L^{-}_{\La}(zq^{-n})^{-1}\Ec_{\La}=\Ec_{\La}\tss L^{-}_{\La}(zq^{-n})^{-1}\Ec_{\La},
\een
and hence
\ben
L_{\La}(z)\tss \Ec_{\La}=\Ec_{\La}\tss L_{\La}(z)\tss\Ec_{\La}.
\een
\ele

\bpf
By the $RLL$-relations \eqref{RLL} we have
\begin{align}
\non
\check{R}_{\La}(z_{1}, \dots, z_{m})\ts L_{1}^{+}&(z_mq^{-2c_m(\La)})\dots L_{m}^{+}(z_1q^{-2c_1(\La)})\\
&=
L_{1}^{+}(z_1q^{-2c_1(\La)})\dots L_{m}^{+}(z_mq^{-2c_m(\La)})\ts\check{R}_{\La}(z_{1}, \dots, z_{m}).
\label{rzll}
\end{align}
Multiply both sides by the element $\check{R}_{0}^{-1}$ from the right and apply consecutive
evaluations by setting $z_1=z,\  z_2=z,\ \dots,\ z_m=z$ to derive from \eqref{fusr} that
\ben
\Ec_{\La}\check{R}_{0}\tss L_{1}^{+}(zq^{-2c_m(\La)})\dots
L_{m}^{+}(zq^{-2c_1(\La)})\check{R}_{0}^{-1}
=L_{1}^{+}(zq^{-2c_1(\La)})\dots L_{m}^{+}(zq^{-2c_m(\La)})\ts\Ec_{\La}.
\een
This shows that the right hand side stays unchanged when it is multiplied by $\Ec_{\La}$
from the left thus proving the first relation.

To prove the second relation, we start
with the counterpart of \eqref{rzll} for the matrix $L^-(z)$ in place of $L^+(z)$.
Then multiply both sides by the inverses of the products of the $L$-matrices
to get
\ben
\bal
\check{R}_{\La}(z_{1}, \dots, z_{m})\ts L_{m}^{-}&(z_1q^{-2c_1(\La)})^{-1}\dots
L_{1}^{-}(z_mq^{-2c_m(\La)})^{-1}\\
&=
L_{m}^{-}(z_mq^{-2c_m(\La)})^{-1}\dots L_{1}^{-}(z_1q^{-2c_1(\La)})^{-1}
\ts\check{R}_{\La}(z_{1}, \dots, z_{m}).
\eal
\een
Now continue with the same argument as for the first relation, but with the consecutive evaluations
$z_1=zq^{-n},\  z_2=zq^{-n},\ \dots,\ z_m=zq^{-n}$ and the application of \eqref{fusr}.
\epf

Introduce an extra copy of the endomorphism algebra $\End\CC^n$ in \eqref{multtpr}
and label it by $0$ to work with the algebra
\beql{multtprze}
\End\CC^n\ot\big(\End\CC^n\big)^{\ot m}\ot\wt\U_q(\wh\gl_n)_{\text{\rm cri}}.
\eeq
The $R$-matrix \eqref{rf} satisfies the Yang--Baxter equation
\ben
R(u/v)R_{02}(w/v)R_{01}(w/u)=R_{01}(w/u)R_{02}(u/v)R(u/v)
\een
which implies
\ben
\check{R}(u/v)R_{02}(w/v)R_{01}(w/u)=R_{02}(w/u)R_{01}(u/v)\check{R}(u/v).
\een
Therefore, by replacing the $L$-matrices with the elements $R_{0i}(z_i)$ in the proof
of Lemma~\ref{lem:LE}, we get its following counterpart.

\ble\label{lem:RE}
We have the relations
\ben
\Ec_{\La}R_{0m}\Big(\frac{zq^{2n+2c_m(\La)}}{w}\Big)\dots
R_{01}\Big(\frac{zq^{2n+2c_1(\La)}}{w}\Big)\Ec_{\La}=
R_{0m}\Big(\frac{zq^{2n+2c_m(\La)}}{w}\Big)\dots R_{01}
\Big(\frac{zq^{2n+2c_1(\La)}}{w}\Big)\Ec_{\La}
\een
and
\ben
\Ec_{\La}R_{01}\Big(\frac{z q^{2c_1(\La)}}{w}\Big)^{-1}\dots
R_{0m}\Big(\frac{z q^{2c_m(\La)}}{w}\Big)^{-1}\Ec_{\La}=
R_{01}\Big(\frac{z q^{2c_1(\La)}}{w}\Big)^{-1}\dots
R_{0m}\Big(\frac{z q^{2c_m(\La)}}{w}\Big)^{-1}\Ec_{\La}.
\een
\ele

Now we will show that
the coefficients of the series $S_{\La}(z)$ belong to the center of the completed quantum
affine algebra at the critical level.

It will be sufficient to verify what $S_{\La}(w)$ commutes with $L^{+}_0(z)$ and $L^-_0(z)$,
regarded as elements of the algebra \eqref{multtprze}. The calculations are similar in both cases
and so we will only consider $L^+_0(z)$ and follow the arguments of \cite[Sec.~3]{fjmr:hs}.
By using \eqref{RLL} and \eqref{RLLpm}
we get
\ben
\bal
L^+_0(z)L_{\La}(w)D_1\dots D_m\Ec_{\La}
&=
R_{01}\Big(\frac{z q^{2c_1(\La)}}{w}\Big)^{-1}\dots
R_{0m}\Big(\frac{z q^{2c_m(\La)}}{w}\Big)^{-1}L_{\La}(w)\\[0.2em]
&\times R_{0m}\Big(\frac{zq^{2n+2c_m(\La)}}{w}\Big)\dots
R_{01}\Big(\frac{zq^{2n+2c_1(\La)}}{w}\Big)
D_1\dots D_m\Ec_{\La}L^+_0(z).
\eal
\een
Therefore, to conclude that $L^+_0(z)\tss S_{\La}(w)=S_{\La}(w)\tss L^+_0(z)$ we need
to show that the trace
\begin{multline}
\tr^{}_{1,\dots,m}\ts
R_{01}\Big(\frac{z q^{2c_1(\La)}}{w}\Big)^{-1}\dots
R_{0m}\Big(\frac{z q^{2c_m(\La)}}{w}\Big)^{-1}L_{\La}(w)\\[0.2em]
\times R_{0m}\Big(\frac{zq^{2n+2c_m(\La)}}{w}\Big)\dots
R_{01}\Big(\frac{zq^{2n+2c_1(\La)}}{w}\Big)
D_1\dots D_m\Ec_{\La}
\label{trace}
\end{multline}
equals $S_{\La}(w)$. Note that the product $D_1\dots D_m$ commutes with both the elements
of the symmetric group $\Sym_m$ and Hecke algebra $\Hc_m$ acting on the tensor product
space $(\CC^n)^{\ot m}$. Hence, using the relation $\Ec_{\La}^2=\Ec_{\La}$ and
applying Lemmas~\ref{lem:LE} and \ref{lem:RE}, we bring \eqref{trace} to the form
\begin{multline}
\tr^{}_{1,\dots,m}\ts
R_{01}\Big(\frac{z q^{2c_1(\La)}}{w}\Big)^{-1}\dots
R_{0m}\Big(\frac{z q^{2c_m(\La)}}{w}\Big)^{-1}
\Ec_{\La}L_{\La}(w)\Ec_{\La}\\[0.2em]
\times R_{0m}\Big(\frac{zq^{2n+2c_m(\La)}}{w}\Big)\dots
R_{01}\Big(\frac{zq^{2n+2c_1(\La)}}{w}\Big)\tss\Ec_{\La}D_1\dots D_m.
\non
\end{multline}
Set
\ben
X= R_{01}\Big(\frac{z q^{2c_1(\La)}}{w}\Big)^{-1}\dots
R_{0m}\Big(\frac{z q^{2c_m(\La)}}{w}\Big)^{-1}
\Ec_{\La}L_{\La}(w)\Ec_{\La}
\een
and
\ben
Y=  R_{0m}\Big(\frac{zq^{2n+2c_m(\La)}}{w}\Big)\dots
R_{01}\Big(\frac{zq^{2n+2c_1(\La)}}{w}\Big)\tss\Ec_{\La}D_1\dots D_m.
\een
Now use the property
\ben
\tr^{}_{1,\dots,m}\ts XY=\tr^{}_{1,\dots,m}\ts X^{t_1\dots t_m}Y^{t_1\dots t_m}.
\een
By using the relation obtained by
applying the transposition $t_1\dots t_m$ to both sides of the second relation in Lemma~\ref{lem:RE},
we get
\begin{multline}
\tr^{}_{1,\dots,m}\ts X^{t_1\dots t_m}Y^{t_1\dots t_m}\\
{}=\tr^{}_{1,\dots,m}\ts\Ec_{\La}^{t_1\dots t_m}L_{\La}^{t_1\dots t_m}(w)
\Bigg(R_{01}\Big(\frac{z q^{2c_1(\La)}}{w}\Big)^{-1}\Bigg)^{t_1}\dots
\Bigg(R_{0m}\Big(\frac{z q^{2c_m(\La)}}{w}\Big)^{-1}\Bigg)^{t_m}\\[0.5em]
\times D_1\dots D_m\tss R_{0m}\Big(\frac{zq^{2n+2c_m(\La)}}{w}\Big)^{t_m}
R_{01}\Big(\frac{zq^{2n+2c_1(\La)}}{w}\Big)^{t_1}.
\non
\end{multline}
As a final step, use the first crossing symmetry relation in \eqref{cs} to conclude that
\ben
\bal
\tr^{}_{1,\dots,m}\ts X^{t_1\dots t_m}Y^{t_1\dots t_m}
{}&=\tr^{}_{1,\dots,m}\ts\Ec_{\La}^{t_1\dots t_m}L_{\La}^{t_1\dots t_m}(w) D_1\dots D_m
\\
{}&=\tr^{}_{1,\dots,m}\ts D_1\dots D_m\Ec_{\La}^{t_1\dots t_m}L_{\La}^{t_1\dots t_m}(w)
\\
{}&=\tr^{}_{1,\dots,m}\ts L_{\La}(w)\Ec_{\La}
D_1\dots D_m
\eal
\een
which coincides with $S_{\La}(w)$. This proves the first part of Theorem~\ref{thm:centr}.

To prove the second part of Theorem~\ref{thm:centr}, it will be enough to verify the relation
\beql{lalapr}
\tr^{}_{1,\dots,m}\ts L_{\La}(z)D_1\dots D_m \Ec_{\La}
=\tr^{}_{1,\dots,m}\ts L_{\La'}(z)D_1\dots D_m \Ec_{\La'},
\eeq
where $\La'=\si_k\La$ is also a standard tableau for $k\in\{1,\dots,m-1\}$.
Set $d_k=d_{k}(\La)$ and note that both operators
$\check{R}_k(q^{-2d_{k}})$ and $\check{R}_k(q^{2d_k})$ are invertible.
By taking into account the property
$P\overline R(x) P=\overline R(x^{-1})^{-1}$
of the $R$-matrix \eqref{rbar} and using \eqref{RLL} we get the relation
\beql{lrch}
L_{\La}(z)\check{R}_k(q^{-2d_k}) = \check{R}_k(q^{-2d_k})L_{\La'}(z).
\eeq
Hence we can write
\begin{align}
L_{\La}(z) D_1\dots D_m \Ec_{\La}
&=  L_{\La}(z)\check{R}_k(q^{-2d_k})\check{R}_k(q^{-2d_k})^{-1}D_1\dots D_m\Ec_{\La}
\non\\
&= \check{R}_k(q^{-2d_k})L_{\La'}(z)\check{R}_k(q^{-2d_k})^{-1}\Ec_{\La} \Ec_{\La}D_1\dots D_m.
\label{prodrld}
\end{align}
Furthermore, the two relations of Lemma~\ref{lem:et} imply
\beql{eq:RE}
\check{R}_k(q^{-2d_{k}})\Ec_{\La'}=\Ec_{\La}\check{R}_k(q^{2d_{k}})
\eeq
so that \eqref{prodrld} takes the form
\begin{multline}
\check{R}_k(q^{-2d_k})L_{\La'}(z)\Ec_{\La'}\check{R}_k(q^{2d_k})^{-1} \Ec_{\La}D_1\dots D_m\\[0.3em]
= \check{R}_k(q^{-2d_k})\Ec_{\La'}L_{\La'}(z)\Ec_{\La'}
\check{R}_k(q^{-2d_k})^{-1} \Ec_{\La}D_1\dots D_m  \\[0.3em]
= \check{R}_k(q^{-2d_k})\Ec_{\La'}L_{\La'}(z)\Ec_{\La'}D_1\dots D_m \check{R}_k(q^{-2d_k})^{-1} \Ec_{\La},
\non
\end{multline}
where we also used Lemma~\ref{lem:LE}.
Now use the cyclic property of trace to write
\ben
\bal
\tr^{}_{1,\dots,m}\ts L_{\La}(z) D_1\dots D_m \Ec_{\La}
&=\tr^{}_{1,\dots,m}\ts\check{R}_k(q^{-2d_k})^{-1} \Ec_{\La}
\check{R}_k(q^{-2d_k})\Ec_{\La'}L_{\La'}(z)\Ec_{\La'}D_1\dots D_m \\
&=\tr^{}_{1,\dots,m}\ts\check{R}_k(q^{-2d_k})^{-1} \Ec_{\La}
\check{R}_k(q^{2d_k})L_{\La'}(z)\Ec_{\La'}D_1\dots D_m \\
&=\tr^{}_{1,\dots,m}\ts\check{R}_k(q^{-2d_k})^{-1}
\check{R}_k(q^{-2d_k})\Ec_{\La'}L_{\La'}(z)\Ec_{\La'}D_1\dots D_m
\eal
\een
which coincides with
\ben
\tr^{}_{1,\dots,m}\ts L_{\La'}(z)\Ec_{\La'}D_1\dots D_m
=\tr^{}_{1,\dots,m}\ts L_{\La'}(z)D_1\dots D_m\Ec_{\La'}
\een
thus proving \eqref{lalapr} and completing the proof of Theorem~\ref{thm:centr}.

\medskip

Since the series $S_{\La}(z)$ does not depend on the standard $\la$-tableau $\La$,
it is unambiguous to set $S_{\la}(z)=S_{\La}(z)$. We will also use the
formula
\beql{sla}
S_{\la}(z)=\tr^{}_{1,\dots,m}\ts T_{\la}(z),\qquad
T_{\la}(z)=\frac{1}{f_{\la}}\sum_{\sh(\La)=\la}\ts L_{\La}(z)D_1\dots D_m\Ec_{\La}.
\eeq

We will now apply Theorem~\ref{thm:centr} to describe
a family of invariants of the
{\em vacuum module at the critical level}
$V_q(\gl_n)$. It is defined as the quotient of $\U_q(\wh\gl_n)_{\text{\rm cri}}$ by the left ideal
generated by all elements $l^-_{ij}[r]$ with $r>0$ and by
the elements $l^-_{i\tss j}[0]-\de_{i\tss j}$ with $i\geqslant j$.
The module $V_q(\gl_n)$ is generated by the vector $\vac$
(the image of $1\in \U_q(\wh\gl_n)_{\text{\rm cri}}$ in the quotient)
such that
$
L^-(u)\tss\vac=I\ts\vac,
$
where $I$ denotes the identity matrix. As a vector space,
$V_q(\gl_n)$ can be identified with the subalgebra
$\Y_q(\gl_n)$ of $\U_q(\wh\gl_n)_{\text{\rm cri}}$ generated by
the coefficients of all series $l^+_{ij}(u)$
subject to the additional relations $l^+_{ii}[0]=1$.
The subspace of invariants of $V_q(\gl_n)$ is defined by
\ben
\z_q(\wh\gl_n)=\{v\in V_q(\gl_n)\ |\ L^-(u)\tss v=I\ts v\};
\een
cf. \cite[Sec.~3.3]{f:lc} and \cite[Sec.~8]{fr:qc}. Some closely related objects were also studied
in \cite{jkmy:cq} and \cite{km:cq} in the context of quantum vertex algebras.
One can regard $\z_q(\wh\gl_n)$ as a subspace of $\Y_q(\gl_n)$.
Moreover, this subspace is closed under the multiplication in
the quantum affine algebra. Therefore, $\z_q(\wh\gl_n)$ can be identified
with a subalgebra of $\Y_q(\gl_n)$.

For any standard $\la$-tableau $\La$ introduce the series $\overline S_{\La}(z)$
with coefficients in $\Y_q(\gl_n)$ by
\ben
\overline S_{\La}(z)=\tr^{}_{1,\dots,m}\ts L^+_1(zq^{-2c_1(\La)})\dots L^+_m(zq^{-2c_{m}(\La)})
D_1\dots D_m \Ec_{\La}.
\een
The following corollary is a generalization of \cite[Cor.~3.3]{fjmr:hs}
and it is verified in the same way.

\bco\label{cor:inv}
The series
$\overline S_{\La}(z)\vac$ does not depend on the standard $\la$-tableau $\La$
and all its coefficients
belong to the algebra of invariants $\z_q(\wh\gl_n)$.
Moreover, the coefficients of all series $\overline S_{\La}(z)$ pairwise commute.
\qed
\eco

\section{Harish-Chandra images}
\label{sec:hch}
\setcounter{equation}{0}

Now we aim to calculate the Harish-Chandra images of the series $S_{\la}(z)$.
Recall the definition of the Harish-Chandra homomorphism as introduced in \cite{fjmr:hs}.
Observing that the coefficients of each generator series \eqref{serlpm} pairwise commute,
introduce a total ordering
on the series such that $l^+_{ij}(u)\prec l^-_{km}(u)$ for all $i,j,k,m$ and
\ben
\bal[]
&l^+_{n\tss 1}(u)\prec l^+_{n-1\ts 1}(u)\prec l^+_{n\ts 2}(u)\prec\dots\prec
l^+_{1\tss 1}(u)\prec\dots\prec
l^+_{n\tss n}(u)\prec l^+_{1\tss 2}(u)\prec\dots\prec l^+_{1\tss n}(u),\\[0.3em]
&l^-_{1\tss n}(u)\prec l^-_{1\ts n-1}(u)\prec l^-_{2\ts n}(u)\prec\dots\prec
l^-_{1\tss 1}(u)\prec\dots\prec
l^-_{n\tss n}(u)\prec l^-_{2\tss 1}(u)\prec\dots\prec l^-_{n\tss 1}(u).
\eal
\een
According to \cite[Prop.~5.1]{fjmr:hs}, the ordered monomials in the generators
form a basis of the
quantum affine algebra $\U_q(\wh\gl_n)$.
Denote by $\U^0$
the subspace of the algebra spanned by
those monomials which do not
contain any generators $l^{\pm}_{ij}[r]$ with $i\ne j$.
Let $x_0$ denote the component
of the linear combination representing the element $x$,
which belongs to $\U^0$.
The mapping $\theta:x\mapsto x_0$ defines the projection
$\theta:\U_q(\wh\gl_n)_{\text{\rm cri}}\to \U^0$. Extending it
by continuity we get the projection $\theta:\wt\U_q(\wh\gl_n)_{\text{\rm cri}}\to \wt\U^0$
to the corresponding completed vector space $\wt\U^0$.

The algebra $\Pi_q(n)$ is defined as
the quotient of the algebra of polynomials in independent variables
$l^+_i[-r]$, $l^-_i[r]$ with $i=1,\dots,n$ and $r=0,1,\dots$
by the relations $l^+_i[0]\tss l^-_i[0]=1$ for all $i$.
The mapping $\eta:\U^0\to \Pi_q(n)$ which takes each ordered monomial
in the generators $l^{\pm}_{ii}[\mp r]$ to the corresponding monomial
in the variables $l^{\pm}_{i}[\mp r]$ by the rule
$l^{\pm}_{ii}[\mp r]\mapsto l^{\pm}_{i}[\mp r]$ extends to an isomorphism
of vector spaces.
Define
the completion $\wt\Pi_q(n)$ of the algebra $\Pi_q(n)$ as
the inverse limit
\ben
\wt\Pi_q(n)=\lim_{\longleftarrow} \Pi_q(n)/I_p, \qquad p>0,
\een
where $I_p$ denotes the ideal of $\Pi_q(n)$ generated by all elements
$l^{-}_{i}[r]$ with $r\geqslant p$; cf. \eqref{compl}. The isomorphism
$\eta$ extends to an isomorphism of the completed
vector spaces $\eta:\wt\U^0\to \wt\Pi_q(n)$.
Thus we get a linear map
\beql{chihom}
\chi:\wt\U_q(\wh\gl_n)_{\text{\rm cri}}\to \wt\Pi_q(n)
\eeq
defined as the composition $\chi=\eta\circ\theta$.
The next proposition provides an analogue of the Harish-Chandra homomorphism
for the quantum affine algebra \cite[Prop.~6.1]{fjmr:hs}.

\bpr\label{prop:hchhom}
The restriction of the map \eqref{chihom} to the center $\Zr_q(\wh\gl_n)$ of the algebra
$\wt\U_q(\wh\gl_n)_{\text{\rm cri}}$ is a homomorphism of commutative algebras
$
\chi:\Zr_q(\wh\gl_n)\to \wt\Pi_q(n).
$
\epr

Combine the generators of the algebra $\Pi_q(n)$ into the series
\ben
l^{+}_{i}(z)=\sum_{r=0}^{\infty}l^{+}_{i}[-r]\tss z^r,\qquad
l^{-}_{i}(z)=\sum_{r=0}^{\infty}l^{-}_{i}[r]\tss z^{-r}
\een
and for $i=1,\dots,n$ set
\ben
x_i(z)=q^{n-2i+1}\ts \frac{l^{+}_{i}(z)\ts l^{-}_{1}(zq^{-n+2})\dots l^{-}_{i-1}(zq^{-n+2i-2})}
{l^{-}_{1}(zq^{-n})\dots l^{-}_{i}(zq^{-n+2i-2})}.
\een
This is a Laurent series in $z$ whose coefficients
are elements of the completed algebra $\wt\Pi_q(n)$.
We are now in a position to state the main result of this section.

\bth\label{thm:hchim}
The image of the series $S_{\la}(z)$ under the Harish-Chandra homomorphism is found by
\ben
\chi:S_{\la}(z)\mapsto \sum_{\sh(\mathcal{T})=\la}
\prod_{\al\in \la}x^{}_{\mathcal{T}(\al)}(zq^{-2c(\al)}),
\een
summed over semistandard tableau $\mathcal{T}$ of shape $\la$ with entries in $\{1,2,\dots, n\}$.
\eth

In the particular case of one-column diagram $\la=(1^m)$, we recover \cite[Thm~6.2]{fjmr:hs}.
The proof in the general case relies on an original argument of \cite{o:qi} which was already
used in the Yangian context in \cite{jkmy:cq}; see also \cite[Ch.~10]{m:so}.
However, the quantum affine algebra case requires some new ingredients with an elaborate use
of the Hecke algebra representations. The proof of Theorem~\ref{thm:hchim} will be given
in the rest of this section.

Given an $m$-tuple $(i)=(i_1,\dots, i_m )$ with
$1\leqslant i_1\leqslant \dots \leqslant i_m\leqslant n$, we let $\al_k$ denote the multiplicity
of the index $k$ in the multiset $\{i_1,\dots, i_m\}$. We thus get the corresponding composition
$\mu_{(i)}=(\al_1,\dots, \al_n)$ of $m$. As with partitions, we will identify $\mu_{(i)}$ with
the associated left-justified diagram of unit boxes such that the top row
contains $\al_1$ boxes, the next row contains $\al_2$ boxes, etc. Similarly,
a $\mu_{(i)}$-tableau is obtained by filling the boxes bijectively with
the numbers $1,\dots,m$. Such a tableau is called {\em row-standard} if the entries
in each row increase from left to right. We let $t^{\mu_{(i)}}$ denote the
row-standard $\mu_{(i)}$-tableau obtained by writing the numbers $1,\dots,m$
in the natural order into the boxes of successive rows starting from the top row.

The Young subgroup $\Sym_{(i)}=\Sym_{\al_1}\times \dots \times \Sym_{\al_n}$ of
the symmetric group $\Sym_m$ preserves
the set of numbers appearing in each row of $t^{\mu_{(i)}}$. The group $\Sym_m$ acts naturally
on the set of $\mu_{(i)}$-tableaux and we introduce the subset $\Dc_{(i)}\subset\Sym_m$ by
\ben
\Dc_{(i)}=\{\om\in \Sym_m\ts|\ts\om\tss t^{\mu_{(i)}} \quad\text{is row-standard}\}.
\een
Clearly $\Dc_{(i)}$ is a set of coset representatives
of $\Sym_{(i)}$ in $\Sym_m$. Given $\si\in\Sym_m$, the unique decomposition $\si=\om\tss\pi$
with $\om\in\Dc_{(i)}$ and $\pi\in\Sym_{(i)}$ is determined by the condition that
the set of elements in each row of $\om\tss t^{\mu_{(i)}}$ coincides with the set
of elements of the corresponding row of $\si\tss t^{\mu_{(i)}}$. Moreover,
by counting the number of inversions in the permutation $\si$ we find that
$l(\si)=l(\om)+l(\pi)$. This implies the relation $T_{\si}=T_{\om}T_{\pi}$
in the Hecke algebra $\Hc_m$.

It will be convenient to use a standard notation for the matrix elements
$A^{i_1\dots i_m}_{j_1\dots j_m}$ of an operator
\ben
A=\sum_{i_1,\dots,i_m,\ts j_1,\dots,j_m}e_{i_1,j_1}\ot \dots \ot e_{i_m,j_m}\ot
A^{i_1\dots i_m}_{j_1\dots j_m}\in(\End\CC^n)^{\ot m}
\een
by writing
\ben
A^{i_1\dots i_m}_{j_1\dots j_m}=\langle i_1,\dots, i_m\tss|\tss A\tss |\tss j_1,\dots,j_m\rangle,
\een
together with the expressions
$
A\tss |\tss j_1,\dots,j_m\rangle$ and $\langle i_1,\dots, i_m\tss|\tss A
$,
interpreted accordingly.

By formula \eqref{rfin} we have the property $\check{R}(e_i\ot e_j)=e_j\ot e_i$
for $i<j$. Therefore,
for an $m$-tuple $(i)=(i_1,\dots, i_m)$ with
$1\leqslant i_1\leqslant \dots \leqslant i_m\leqslant n$ and any $\om\in \Dc_{(i)}$
we have
\ben
\check{R}_{\om}\tss|\tss i_1,\dots, i_m\rangle=|\tss i_{\om(1)},\dots, i_{\om(m)}\rangle
\Fand
\langle i_1,\dots, i_m\tss |\tss\check{R}_{\om^{-1}}=\langle i_{\om(1)},\dots, i_{\om(m)}\tss|.
\een
Formula \eqref{sla} for $S_{\la}(z)$ then implies
\begin{multline}
S_{\la}(z)
=\sum_{j_1,\dots, j_m=1}^{n}\langle j_1,\dots, j_m\tss|\tss T_{\la}(z)\tss |\tss j_1,\dots,j_m\rangle\\
=\sum_{i_1\leqslant\dots\leqslant i_m}
\sum_{\om\in \Dc_{(i)}}\langle i_1,\dots, i_m\tss|\tss
\check{R}_{\om^{-1}}T_{\la}(z)\check{R}_{\om} \tss|\tss i_1,\dots, i_m\rangle.
\non
\end{multline}

\ble\label{lem:RT}
For any $\si\in \Sym_m$ and $\la \vdash m$ we have
\ben
\check{R}_{\si}{T}_{\la}(z)={T}_{\la}(z)\check{R}_{\si}.
\een
\ele

\bpf
We will use the isomorphism \eqref{isomhecke}
between the algebra $\Hc_m$ and the direct sum of matrix algebras.
We will show that the element ${T}_{\la}(z)$ commutes with
the image $\Ec^{\mu}_{S T}$ of an arbitrary matrix unit $e^{\mu}_{S T}$
under the action of the Hecke algebra $\Hc_m$ in $(\CC^{n})^{\otimes m},$
where $S$ and $T$ are standard $\mu$-tableaux. By Lemma~\ref{lem:LE} we have
\ben
{T}_{\la}(z)=\frac{1}{f_{\la}}\sum_{\operatorname{sh}(\La)=\la}
\Ec^{\la}_{\La}D_1\dots D_mL_{\La}(z)\Ec^{\la}_{\La}.
\een
Hence, if $\mu \neq \la,$ then
\ben
\Ec^{\mu}_{S T} {T}_{\la}(z)={T}_{\la}(z)\tss \Ec^{\mu}_{S T}=0.
\een
Now suppose $\mu =\la$. In the case where $S=T$ we have
\ben
\Ec^{\la}_{S}{T}_{\la}(z)={T}_{\la}(z) \Ec^{\la}_{S}=
\frac{1}{f_{\la}}\ts\Ec^{\la}_{S}D_1\dots D_mL_{S}(z)\Ec^{\la}_{S}.
\een
Finally, suppose $S$ and $T$ are different standard tableaux of shape $\la$.
It suffices to consider the case where $T$ is obtained from $S$
by swapping the entries $k$ and $k+1$ for some $k \in\{1, \dots, m-1\}$ so that
$T=\si_{k} S$. In this case, Lemma \ref{lem:et} implies
\ben
{T}_{\la}(z) \Ec^{\la}_{S T}
=\frac{1}{\de}\ts{T}_{\la}(z)\Ec_{S}^{\la}\check{R}_k(q^{2d_k})
=\frac{1}{\de\tss f_{\la}}\ts L_{S}(z)D_1\dots D_m\Ec^{\la}_{S}\check{R}_k(q^{2d_k}),
\een
where we set
\ben
d_k=d_k(S)\Fand\de=\sqrt{1-\frac{1}{[d_k]_q^2}}\ .
\een
Due to \eqref{eq:RE},
we have
\ben
{T}_{\la}(z) \Ec^{\la}_{S T}
=\frac{1}{\de\tss f_{\la}}\ts L_{S}(z)D_1\dots D_m\check{R}_k(q^{-2d_k})\Ec^{\la}_{T}
=\frac{1}{\de\tss f_{\la}}\ts L_{S}(z)\check{R}_k(q^{-2d_k})\Ec^{\la}_{T}D_1\dots D_m.
\een
Now use \eqref{lrch} to write this expression as
\ben
\frac{1}{\de\tss f_{\la}}\ts\check{R}_k(q^{-2d_k})L_{T}(z)\Ec^{\la}_{T}D_1\dots D_m
=\frac{1}{\de\tss f_{\la}}\ts\check{R}_k(q^{-2d_k})\Ec^{\la}_{T}L_{T}(z)\Ec^{\la}_{T}D_1\dots D_m
\een
which equals
\ben
\frac{1}{\de}\ts\check{R}_k(q^{-2d_k})\Ec^{\la}_{T}{T}_{\la}(z)
=\frac{1}{\de}\ts\Ec^{\la}_{S}\check{R}_k(q^{2d_k}{T}_{\la}(z)=\Ec^{\la}_{S T}{T}_{\la}(z),
\een
thus completing the proof.
\epf

Lemma~\ref{lem:RT} allows us to write the expression for
$S_{\la}(z)$ in the form
\ben
\bal
S_{\la}(z)&=\sum_{i_1\leqslant\dots\leqslant i_m}
\sum_{\om\in \Dc_{(i)}}\langle i_1,\dots, i_m\tss|\tss T_{\la}(z)
\check{R}_{\om^{-1}}\check{R}_{\om} \tss|\tss i_1,\dots, i_m\rangle\\
&=\frac{1}{f_{\la}}\sum_{i_1\leqslant\dots\leqslant i_m}
\sum_{\om\in \Dc_{(i)}}\sum_{\sh(\La)=\la}
\langle i_1,\dots, i_m\tss|\tss L_{\La}(z)D_1\dots D_m
\Ec_{\La}\check{R}_{\om^{-1}}\check{R}_{\om}\tss|\tss i_1,\dots, i_m\rangle.
\eal
\een
By formula \eqref{tauchi},
\ben
\Ec_{\La}\check{R}_{\om^{-1}}\check{R}_{\om}
=\sum\limits_{\si\in \Sym_m}\tau\left(e_{\La}T_{\om^{-1}}T_{\om}T_{\si^{-1}}\right)\check{R}_{\si}
=\sum\limits_{\si\in \Sym_m}\frac{1}{c_{\la}}
\chi_{\la}(e_{\La}T_{\om^{-1}}T_{\om}T_{\si^{-1}})\check{R}_{\si}.
\een
Therefore, $S_{\la}(z)$ equals
\ben
\frac{1}{f_{\la}{c_{\la}}}\sum_{i_1\leqslant\dots\leqslant i_m}
\sum_{\substack{\om\in \Dc_{(i)}\\\si\in \Sym_m}}\sum_{\sh(\La)=\la}
\chi_{\la}\left(e_{\La}T_{\om^{-1}}T_{\om}T_{\si^{-1}}\right)
\langle i_1,\dots, i_m\tss|\tss L_{\La}(z)D_1\dots D_m
\check{R}_{\si}\tss|\tss i_1,\dots, i_m\rangle.
\een
Using the decomposition $\Sym_m=\Dc_{(i)}\Sym_{(i)}$, we can write this expression
in the form
\begin{multline}
\frac{1}{f_{\la}{c_{\la}}}\sum_{i_1\leqslant\dots\leqslant i_m}
\sum_{\substack{\om\in \Dc_{(i)}}}
\sum_{\substack{\om'\in \Dc_{(i)}\\\pi\in \Sym_{(i)}}}\sum_{\sh(\La)=\la}
\chi_{\la}\big(e_{\La}T_{\om^{-1}}T_{\om}T_{\pi^{-1}}T_{\om'^{-1}}\big)\\[-1em]
{}\times\langle i_1,\dots, i_m\tss|\tss L_{\La}(z)D_1\dots D_m
\check{R}_{\om'}\check{R}_{\pi}\tss|\tss i_{1},\dots, i_{m}\rangle.
\non
\end{multline}
For an $m$-tuple $(i_1,\dots, i_m)$ with $1\leqslant i_1\leqslant \dots \leqslant i_m \leqslant n$
and $\pi\in \Sym_{(i)}$ one easily verifies that
\ben
\check{R}_{\pi}\tss|\tss i_{1},\dots, i_{m}\rangle
=q^{l(\pi)}\tss|\tss i_{\pi(1)},\dots, i_{\pi(m)}\rangle
=q^{l(\pi)}\tss|\tss i_{1},\dots, i_{m}\rangle.
\een
Moreover, we have
\ben
\check{R}_{\om'}\check{R}_{\pi}\tss|\tss i_{1},\dots, i_{m}\rangle
=q^{l(\pi)}\check{R}_{\om'}\tss|\tss i_{1},\dots, i_{m}\rangle
=q^{l(\pi)}\tss|\tss i_{\om'(1)},\dots, i_{\om'(m)}\rangle.
\een
Hence
\begin{multline}
S_{\la}(z)=\frac{1}{f_{\la}{c_{\la}}}\sum_{i_1\leqslant\dots\leqslant i_m}
\sum_{\substack{\om\in \Dc_{(i)}}}
\sum_{\substack{\om'\in \Dc_{(i)}\\\pi\in \Sym_{(i)}}}\sum_{\sh(\La)=\la}
q^{l(\pi)}\chi_{\la}\big(e_{\La}T_{\om^{-1}}T_{\om}T_{\pi^{-1}}T_{\om'^{-1}}\big)\\[-1em]
{}\times\langle i_1,\dots, i_m\tss|\tss L_{\La}(z)D_1\dots D_m
\tss|\tss i_{\om'(1)},\dots, i_{\om'(m)}\rangle
\non
\end{multline}
which equals
\begin{multline}
\frac{1}{f_{\la}c_{\la}}\sum_{i_1\leqslant\dots\leqslant i_m}
\sum_{\substack{\om\in \Dc_{(i)}}}
\sum_{\substack{\om'\in \Dc_{(i)}\\\pi\in \Sym_{(i)}}}
\sum_{\sh(\La)=\la}\sum_{j_1,\dots, j_m}q^{l(\pi)}\chi_{\la}\big(e_{\La}
T_{\om^{-1}}T_{\om}T_{\pi^{-1}}T_{\om'^{-1}}\big)\\
{}\times l_{i_{1},j_1}^{+}(zq^{-2c_{1}(\La)})\dots l_{i_{m},j_m}^{+}(zq^{-2c_{m}(\La)})
\ts \wt{l}_{j_m ,i_{\om'(m)}}(zq^{-n-2c_{m}(\La)})
 \dots \wt{l}_{j_1 ,i_{\om'(1)}}(zq^{-n-2c_{1}(\La)}),
\non
\end{multline}
where we used the entries of the matrix $L^-(z)^{-1}D=[\tss\wt l_{ij}(z)]$.
As with the particular case of one-column diagram $\la$ considered in \cite{fjmr:hs},
the definition of the homomorphism \eqref{chihom} implies that
a nonzero contribution to the Harish-Chandra
image of $S_{\la}(z)$ comes only from the terms with $\om'=1$ and $j_k=i_k$ for all $k=1,\dots,m$.
Therefore, the image is given by
\begin{multline}
\frac{1}{f_{\la}c_{\la}}\sum_{i_1\leqslant\dots\leqslant i_m}
\sum_{\substack{\om\in \Dc_{(i)}\\\pi\in \Sym_{(i)}}}\sum_{\sh(\La)=\la}
{ q^{l(\pi)}}\chi_{\la}\big(e_{\La}T_{\om^{-1}}T_{\om}T_{\pi^{-1}}\big)\\
{}\times l_{i_1,i_1}^{+}(zq^{-2c_{1}(\La)})\dots l_{i_{m},i_m}^{+}(zq^{-2c_{m}(\La)})
\ts\wt{l}_{i_m ,i_{m}}(zq^{-n-2c_{m}(\La)})\dots
 \wt{l}_{i_1 ,i_{1}}(zq^{-n-2c_{1}(\La)})
\non
\end{multline}
which coincides with
\ben
\frac{1}{f_{\la}c_{\la}}\sum_{i_1\leqslant\dots\leqslant i_m}
\sum_{\substack{\om\in \Dc_{(i)}\\\pi\in \Sym_{(i)}}}\sum_{\sh(\La)=\la}
{ q^{l(\pi)}}\chi_{\la}\big(e_{\La}T_{\om^{-1}}T_{\om}T_{\pi}\big)
x_{i_1}(zq^{-2c_{1}(\La)})\dots x_{i_m}(zq^{-2c_{m}(\La)}).
\een

Given an $m$-tuple $(i_1,\dots, i_m)$ with $i_1\leqslant i_2\leqslant \dots \leqslant i_m$,
we let $\Tc=i(\La)$ denote the
tableau obtained from a standard $\la$-tableau $\La$ by replacing
the entry $r$ with $i_r$ for $r=1,\dots, m$.
The entries of $\Tc$ weakly increase
along the rows and down the columns.
Changing the order of summation in the above expression
for the Harish-Chandra image of $S_{\la}(z)$, we can write it as
\ben
\bal
\frac{1}{f_{\la}{c_{\la}}}\sum_{i_1\leqslant\dots\leqslant i_m}\ts
\sum_{\Tc,\ts\sh(\Tc)=\la}\Bigg(
\sum_{\substack{\sh(\La)=\la\\i(\La)=\Tc}}\ts
\sum\limits_{\substack{\om\in \Dc_{(i)}\\\pi\in \Sym_{(i)}}}{ q^{l(\pi)}}
\chi_{\la}\big(e_{\La}T_{\om^{-1}}T_{\om}T_{\pi}\big)
\Bigg)\prod_{\al\in \la}x^{}_{\Tc(\al)}(zq^{-2c(\al)}),
\eal
\een
where $\Tc$ runs over $\la$-tableaux with entries in $\{1,\dots,n\}$ such that
the entries of $\Tc$ weakly increase
along the rows and down the columns.
To complete the proof of Theorem~\ref{thm:hchim}, it will therefore be enough to show that
for a given $\Tc$ we have
\beql{idenchi}
\sum_{\substack{\sh(\La)=\la\\i(\La)=\Tc}}\ts
\sum\limits_{\substack{\om\in \Dc_{(i)}\\\pi\in \Sym_{(i)}}}{ q^{l(\pi)}}
\chi_{\la}\big(e_{\La}T_{\om^{-1}}T_{\om}T_{\pi}\big)
=\begin{cases}f_{\la}\tss c_{\la}\quad&\text{if $\Tc$ is semistandard},\\
     0\quad &\text{otherwise.}
   \end{cases}
\eeq
Let $\Hc_{(i)}$ denote the subalgebra of $\Hc_m$ corresponding to
the Young subgroup $\Sym_{(i)}$ of $\Sym_m$.
Introduce the notation
\ben
s_{(i)}=\sum\limits_{\pi\in \Sym_{(i)}}q^{l(\pi)}T_{\pi}\in\Hc_{(i)}\Fand
e^{}_{\Tc}=\sum\limits_{\substack{\sh(\La)=\la\\i(\La)=\Tc}}e^{}_{\La}.
\een
Then $\CC s_{(i)}$ is the trivial representation of $\Hc_{(i)}$, and the induced representation
$\Hc_{m}\ot_{\Hc_{(i)}}\CC s_{(i)}$ is the permutation representation of $\Hc_m$
corresponding to the composition $\mu_{(i)}=(\al_1,\dots, \al_n)$ of $m$.
The element $s_{(i)}$ is proportional to a central idempotent so that
\beql{misq}
s_{(i)}^2=\prod_{r=1}^n{[\al_r]_q!\ts q^{\al_r(\al_r-1)/2}}\ts s_{(i)}.
\eeq

The left hand side of \eqref{idenchi} can now be written as
\beql{chico}
\sum\limits_{\om\in \Dc_{(i)}} \chi_{\la}\big(e_{\Tc}T_{\om^{-1}}T_{\om}s_{(i)}\big)=
\sum\limits_{\om\in \Dc_{(i)}} \chi_{\la}\big(T_{\om}s_{(i)}e_{\Tc}T_{\om^{-1}}\big).
\eeq
We will use the homomorphism $\vp_{\la}$ associated with the $\Hc_m$-module $V_{\la}$; see \eqref{phila}.

\ble\label{lem:ASchur}
We have the relation
\beql{phichi}
\sum_{\om\in \Dc_{(i)}}\vp_{\la}\big(T_{\om}s_{(i)}e_{\Tc}T_{\om^{-1}}\big)
=\prod_{r=1}^n\frac{1}{[\al_r]_q!\ts q^{\al_r(\al_r-1)/2}}\ts
c_{\la}\ts\chi_{\la}(s_{(i)}e_{\Tc})\ts{\rm id}^{}_{V_{\la}}.
\eeq
\ele

\bpf
Denote by
$V_{\Tc}$ the linear span of the basis vectors $v_{\La}\in V_{\la}$ such that
$i(\La)=\Tc$. It follows from \eqref{eq:Young basis} that the subspace $V_{\Tc}$
of $V_{\la}$ is invariant under
the action of the subalgebra $\Hc_{(i)}$.
Therefore, if a standard $\la$-tableau $\La$ is such that $i(\La)=\Tc$, then
$\vp_{\la}(e_{\Tc}s_{(i)})(v_{\La})=\vp_{\la}(s_{(i)})(v_{\La})$.
Otherwise, if $i(\La)\neq \Tc$ then
\ben
\bal
\vp_{\la}(e_{\Tc}s_{(i)})(v_{\La})&=\sum_{\sh(\Ga)=\la}
\langle s_{(i)}v_{\La},v^{}_{\Ga}\rangle e^{}_{\Tc}v^{}_{\Ga}
=\sum_{\sh(\Ga)=\la}\langle v_{\La},s_{(i)}^{*}v^{}_{\Ga}\rangle
e^{}_{\Tc}v^{}_{\Ga}\\
&=\sum\limits_{\substack{\sh(\Ga)=\la\\i(\Ga)=\Tc}}
\langle v_{\La},s_{(i)}v^{}_{\Ga}\rangle e^{}_{\Tc}v^{}_{\Ga}
+\sum\limits_{\substack{\sh(\Ga)=\la\\i(\Ga)\ne\Tc}}
\langle v_{\La},s_{(i)}v^{}_{\Ga}\rangle e^{}_{\Tc}v^{}_{\Ga}=0.
\eal
\een
On the other hand, we also have
\ben
\vp_{\la}(s_{(i)}e_{\Tc})(v_{\La})=\begin{cases}
   0\qquad&\text{if\ \  $i(\La)\neq \Tc$},\\
    \vp_{\la}( s_{(i)})(v_{\La})\qquad&\text{if\ \  $i(\La)= \Tc$.}
   \end{cases}
\een
Therefore, we may conclude that $\vp_{\la}(s_{(i)}e_{\Tc})=\vp_{\la}(e_{\Tc}s_{(i)})$.

Furthermore, taking into account \eqref{misq}, we can write
\ben
\sum_{\om\in \Dc_{(i)}}\vp_{\la}(T_{\om}s_{(i)}e_{\Tc}T_{\om^{-1}})
=\prod\limits_{r=1}^n\frac{1}{[\al_r]_q!q^{\al_r(\al_r-1)/2}}
\sum_{\om\in \Dc_{(i)}}\vp_{\la}\big(T_{\om}s_{(i)}e_{\Tc}s_{(i)}T_{\om^{-1}}\big).
\een
Since $T_{\pi}s_{(i)}=q^{l(\pi)}s_{(i)}$ for $\pi\in \Sym_{(i)}$, we have
\begin{multline}
\sum_{\om\in \Dc_{(i)}}\vp_{\la}\big(T_{\om}s_{(i)}e_{\Tc}s_{(i)}T_{\om^{-1}}\big)
=\sum\limits_{\substack{\om\in \Dc_{(i)}\\\pi\in \Sym_{(i)}}}
\vp_{\la}\big(q^{l(\pi)}T_{\om}s_{(i)}e_{\Tc}T_{\pi^{-1}}T_{\om^{-1}}\big)\\[-1em]
{}=\sum\limits_{\substack{\om\in \Dc_{(i)}\\\pi\in \Sym_{(i)}}}
\vp_{\la}\big(T_{\om}T_{\pi}s_{(i)}e_{\Tc}T_{\pi^{-1}}T_{\om^{-1}}\big).
\non
\end{multline}
Now observe that
$\{T_{\om}T_{\pi}\tss|\tss\om\in \Dc_{(i)}, \pi\in \Sym_{(i)}\}$
and $\{T_{\pi^{-1}}T_{\om^{-1}}\tss|\tss\om\in \Dc_{(i)}, \pi\in \Sym_{(i)}\}$
are dual bases of $\Hc_m$.
Hence \eqref{iuschur} yields
\ben
\sum_{\om\in \Dc_{(i)}}\vp_{\la}\big(T_{\om}s_{(i)}e_{\Tc}T_{\om^{-1}}\big)
=\prod\limits_{r=1}^n\frac{1}{[\al_r]_q!q^{\al_r(\al_r-1)/2}}\ts
c_{\la}\ts\tr^{}_{V_{\la}}\big(\vp_{\la}(s_{(i)}e_{\Tc})\big)
{\rm id}^{}_{V_{\la}},
\een
which equals the right hand side of \eqref{phichi}.
\epf

Returning to the proof of \eqref{idenchi}, note that by Lemma~\ref{lem:ASchur} the expression
\eqref{chico} equals
\ben
\sum\limits_{\om\in \Dc_{(i)}} \chi_{\la}
\left(T_{\om}s_{(i)}e_{\Tc}T_{\om^{-1}}\right)=
\prod_{r=1}^n\frac{1}{[\al_r]_q!\ts q^{\al_r(\al_r-1)/2}}\ts c_{\la}\tss f_{\la}\ts
\chi_{\la}(s_{(i)}e_{\Tc}).
\een
Thus, the following lemma will complete the proof of Theorem~\ref{thm:hchim}.

\ble\label{lem:skew}
We have
\ben
\chi_{\la}(s_{(i)}e_{\Tc})
=\begin{cases} \prod\limits_{r=1}^n[\al_r]_q!\ts q^{\al_r(\al_r-1)/2}
\quad&\text{if $\Tc$ is semistandard},\\
     0\quad &\text{otherwise.}
     \end{cases}
\een
\ele

\bpf
Write
\ben
s_{(i)}e_{\Tc}=
\sum\limits_{\substack{\pi\in \Sym_{(i)}}}{ q^{l(\pi)}}T_{\pi}e_{\Tc}
=
\prod_{r=1}^n\sum\limits_{\substack{\si_r\in \Sym_{\al_r}}} q^{l(\si_r)}T_{\si_r}e_{\Tc}.
\een
Hence
\ben
\chi_{\la}(s_{(i)}e_{\Tc})
=\sum\limits_{\substack{\sh(\La)=\la\\i(\La)=\Tc}}
\prod_{r=1}^n\sum\limits_{\substack{\si_r\in \Sym_{\al_r}}} q^{l(\si_r)}
\langle T_{\si_r}v_{\La},v_{\La}\rangle=\prod _{r=1}^n\sum_{\si_r\in \Sym_{\al_{r}}}q^{l(\si_r)}
\chi_{\om_r}(T_{\si_{r}}),
\een
where $\om_r$ is the skew diagram which
consists of the boxes of $\Tc$ occupied by $r$, and
$\chi_{\om_r}$ denotes the skew character of $\Hc_{\al_r}$ associated with $\om_r$.
Since $q^{l(\si_r)}=\chi_{\iota}(T_{\si_{r}^{-1}})$ for the trivial
representation $\iota$ of the Hecke algebra $\Hc_{\al_r}$, we can write
\ben
\chi_{\la}(s_{(i)}e_{\Tc})=\prod _{r=1}^n\sum_{\si_r\in \Sym_{\al_{r}}}\chi_{\iota}(T_{\si_{r}^{-1}})
\chi_{\om_r}(T_{\si_{r}}).
\een
The multiplicity of the trivial representation
$\iota$ in the skew representation of $\Hc_{\alpha_{r}}$
associated with $\om_{r} $ is zero unless $\om_{r}$
does not contain two boxes in the same column, in which case the multiplicity is $1$.
Then by \eqref{orthcha},
\ben
\sum_{\si_r\in \Sym_{\al_{r}}}q^{l(\si_r)}
\chi_{\om_r}(T_{\si_{r}})
=\sum_{\si_r\in \Sym_{\al_{r}}}
\chi_{\iota}(T_{\si_{r}})\chi_{\iota}(T_{\si_{r}^{-1}})
=[\al_r]_q!\ts q^{\al_r(\al_r-1)/2}.
\een
This proves the lemma and completes the proof of Theorem~\ref{thm:hchim}.
\epf

\bre\label{rem:qchar}
By replacing each series $x^{}_{i}(zq^{-2c})$ with a formal variable $x^{}_{i,c}$ we find
that the Harish-Chandra image of $S_{\la}(z)$ found in Theorem~\ref{thm:hchim}
takes the form of the $q$-character of the
evaluation module $L(\la^{\circ})$ over $\U_q(\wh\gl_n)$ with
$\la^{\circ}=(\la_1,\dots,\la_m,0,\dots,0)$; cf. \cite[Sec.~7.4]{bk:rs}
and \cite[Sec.~4.5]{fm:ha}. This agrees with the observation
that the subspace $\Ec^{\la}_{\La}(\CC^n)^{\ot m}$ carries the irreducible representation
of $\U_q(\gl_n)$ isomorphic to $L(\la^{\circ})$.
\qed
\ere

Consider the $q$-deformed Wakimoto modules
over the quantum affine algebras in type $A$ constructed in \cite{aos:fb}.
At the critical level, such modules $\Wc(\vk(z))$ are parameterized
by power series $\vk^+_1(z),\dots,\vk^+_n(z)$ in $z$
and a power series $\vk^-(z)$ in $z^{-1}$; see \cite[Sec.~7]{fjmr:hs}
for the definition in terms of the $RLL$ presentation.
It was shown in {\em loc. cit.} that, under a suitable
identification of parameters, the Harish-Chandra images of the elements
of the center of the algebra $\wt\U_q(\wh\gl_n)$
at the critical level coincide
with the eigenvalues of the central elements acting in $\Wc(\vk(z))$. Namely, the generator series
of the algebra $\Pi_q(n)$
should be specialized as
\ben
l^{+}_{i}(z)\mapsto\vk^+_i(z)\Fand l^{-}_{i}(z)\mapsto\vk^-(z)
\een
for $i=1,\dots,n$. The following corollary is therefore immediate from
Theorem~\ref{thm:hchim}.

\bco\label{cor:waki}
The eigenvalues of the coefficients of the series
$S_{\la}(z)$ in the module $\Wc(\vk(z))$ are found by
\ben
S_{\la}(z)\mapsto\sum_{\sh(\mathcal{T})=\la}
\prod_{\al\in \la}\vk^{}_{\mathcal{T}(\al)}(zq^{-2c(\al)}),
\een
where
\ben
\vk_i(z)=q^{n-2i+1}\ts \vk^+_{i}(z)\vk^-(zq^{-n})^{-1},\qquad i=1,\dots,n.
\een
\eco

\newpage

\small

\noindent
N.J.:\newline
Department of Mathematics\\
North Carolina State University, Raleigh, NC 27695, USA\\
jing@math.ncsu.edu

\vspace{5 mm}

\noindent
M.L.:\newline
School of Mathematical Sciences\\
South China University of Technology\\
Guangzhou, Guangdong 510640, China\\
mamliu@scut.edu.cn

\vspace{5 mm}

\noindent
A.M.:\newline
School of Mathematics and Statistics\newline
University of Sydney,
NSW 2006, Australia\newline
alexander.molev@sydney.edu.au

\end{document}